 \documentclass[SecEq,CM]{degruyter-crelle} 

\title[$L$-func. of cert. exp. sums]{On $L$-functions of certain exponential sums}

\author{Jun}{Zhang}{}{Tianjin}
\author{Weiduan}{Feng}{}{Shanghai}

\contact{Chern Institute of Mathematics, Nankai University, Tianjin 300071, P.R. China}{zhangjun04@mail.nankai.edu.cn}
\contact{Department of Mathematics, Shanghai Jiao Tong University, Shanghai 200240, P.R. China}{d.d.feng@163.com}

\researchsupported{The research of Jun Zhang is supported by the National Key Basic Research Program of China (Grant No. 2013CB834204), and the National Natural Science Foundation of China (Nos. 61171082, 10990011 and 60872025).}

\usepackage{epsfig}
\usepackage{color}
\usepackage{amsfonts}
\usepackage{amssymb}
\usepackage{amsmath}
\usepackage{graphicx}

\theoremstyle{plain}
 \newtheorem{thm}{Theorem}[section]
\newtheorem{lem}{Lemma}[section]
\newtheorem{prop}{Proposition}[section]
\newtheorem{cor}{Corollary}[section]

\theoremstyle{definition}
 \newtheorem{rem}{Remark}[section]

\newtheorem*{pf}{Proof}
\newtheorem{defi}{Definition}[section]
\numberwithin{equation}{section}

\def\f#1{{\mathbb{F}}_{#1}}
\def\q#1{{\mathbb{Q}}_{#1}}
\def\z#1{{\mathbb{Z}}_{#1}}
\begin{document}

\begin{abstract}
Let $\f{q}$ denote the finite field of order $q$ (a power of a prime $p$). We study the $p$-adic valuations for zeros  of $L$-functions associated with exponential sums of the following family of
Laurent polynomials
\begin{equation*}
f(x_1,x_2,\cdots,x_{n+1})=a_1x_{n+1}(x_1+{1\over x_1})+\cdots+a_{n}x_{n+1}(x_{n}+{1\over x_{n}})+a_{n+1}x_{n+1}+{1\over x_{n+1}}
\end{equation*}
where $a_i\in \f{q}^*,\,i=1,2,\cdots,n+1$. When $n=2$, the estimate of the associated exponential sum appears in Iwaniec's work,  and Adolphson and Sperber gave complex absolute values for zeros of the corresponding $L$-function.
 Using the decomposition theory of Wan,
we determine the generic Newton polygon ($q$-adic values of the reciprocal zeros) of the $L$-function.
Working on the chain level version of Dwork's trace formula and using Wan's decomposition theory, we are able to give an  explicit Hasse polynomial
for the generic Newton polygon in low dimensions, i.e., $n\leq 3$.

\end{abstract}


\section{Introduction}
$L$-functions have been a powerful tool to investigate exponential sums in number theory. To estimate an exponential sum, people are interested in the zeros and poles of the corresponding $L$-function. Mathematicians study the number of these zeros and poles, the complex absolute values, $l$-adic absolute values of them for primes $l\neq p$, and $p$-adic absolute values of them, especially for some interesting varieties and exponential sums~\cite{FW1,FW2,RW1,Spe1,Spe2}. Deligne's theorem gives the general information for complex absolute values of the zeros and poles of $L$-function. For $l$-adic absolute values, it is well-known that if $l\neq p$ then all the zeros and poles have $l$-adic absolute value $1$. However, for $p$-adic absolute values, it is still very mysterious~\cite{Phong,Maz,Mil,Mum,Zhu04,Zhu12}, especially in higher dimensions.

In this paper, we consider the following family of Laurent polynomials
\begin{equation}\label{poly}
f(x_1,x_2,\cdots,x_{n+1})=a_1x_{n+1}(x_1+{1\over x_1})+\cdots+a_{n}x_{n+1}(x_{n}+{1\over x_{n}})+a_{n+1}x_{n+1}+{1\over x_{n+1}}
\end{equation}
where $a_i\in \f{q}^*,\,i=1,2,\cdots,n+1$.  The exponential sum associated to $f$ is defined to be
$$
S_k^*(f)=\sum_{x_1,\ldots,x_n\in
\f{q^k}^*}\zeta_p^{{\rm{Tr}}_kf(x_1,\ldots,x_n)},
$$
where $\zeta_p$ is a fixed primitive
$p$-th root of unity in the complex numbers
and ${\rm{Tr}}_k$ denotes the trace map from the $k$-th extended field $\f{q^k}$ to the prime
field $\f{p}$.  When $n=2$, the estimate of exponential sum $S_k^*(f)$
is vital in analytic number theory. It appears in Iwaniec's work on
small eigenvalues of the Laplace-Beltrami operator acting on
automorphic functions with respect to the group $\Gamma_0(p)$.

To understand the sequence $S_k^*(f)\in
\mathbb{Q}(\zeta_p)$ $(1\le k<\infty)$ of algebraic integers, we
study the $L$-function associated to $S_k^*(f)$
$$
L^*(f,T)=\exp\left(\sum_{k=1}^{\infty} S_k^*(f){T^k\over k}\right).
$$

By the theorem of Adolphson and Sperber \cite{AS89}, the $L$-function $L^*(f,T)^{(-1)^n}$ for non-degenerate $f$ is a polynomial of degree $(n+1)!\mathrm{Vol}(\Delta)$, where $\Delta=\Delta(f)$
is the Newton polyhedron of $f$ defined explicitly later.
As the origin is an interior point of the Newton polyhedron $\Delta$, by the theorem of Denef and Loeser~\cite{DL}, we have
\begin{thm}
Assume that $f$ is non-degenerate, that is,
$$
\prod_{(c_1,c_2,\cdots,c_n)\in \{\pm 1\}^n}\left(2c_1a_1+2c_2a_2+\cdots+2c_na_n+a_{n+1}\right)\ne 0.
$$
Then, the $L$-function $L^*(f,T)$ associated to the exponential sum $S_k^*(f)$ is pure of weight $n+1$, i.e.,
\[
  L^*(f,T)^{(-1)^n}=\prod_{i=1}^{(n+1)!\mathrm{Vol}(\Delta(f))}(1-\alpha_{i}T)
\]
with the complex absolute value $|\alpha_{i}|=q^{(n+1)/2}$.
\end{thm}
On the other hand, for each $l$-adic absolute value
$|\cdot|_l$ with prime $l\ne p$, the reciprocal zeros $\alpha_i$ are
$l$-adic units: $|\alpha_i|_l=1.$ So we next study the $p$-adic slopes of the reciprocal zeros of $L^*(f,T)^{(-1)^{n}}$ of such a family of exponential sums for the remaining prime $p$.
When $n=1$, such non-degenerate Laurent polynomials $f$ are always ordinary and the Hodge polygon is very clear. So, we only consider $n\ge 2$. One of our main results in this paper is

\begin{thm}\label{main} Assume that $f$ of the form (\ref{poly}) is non-degenerate.

(i) The polynomial $L^*(f,T)^{(-1)^n}$ has degree $2^{n+1}$.

(ii) There exists a non-zero polynomial $h_p(\Delta)(a_1,a_2,\cdots,a_{n+1})\in \f{p}[a_1,a_2,\cdots,a_{n+1}]$ such that if $f$ has coefficients $a_1,a_2,\cdots,a_{n+1}\in \f{q}^*$ verifying $h_p(\Delta)(a_1,a_2,\cdots,a_{n+1})\neq 0$, then for each $k=0,1,\cdots,n+1$, the number of reciprocal zeros of $L^*(f,T)^{(-1)^n}$ with $q$-adic slope $k$ is
\[
 \binom{n+1}{k},
\]
and for any rational number $k \not\in \{0, 1, ..., n+1\}$, there is no reciprocal zero of $L^*(f,T)^{(-1)^n}$ with $q$-adic slope $k$.
\end{thm}

We identify a Laurent polynomial $f$ of the form~(\ref{poly}) with the vector of its coefficients $a=(a_1,a_2,\cdots,a_{n+1})$, as they are one-to-one correspondent to each other.
Let $\mathcal {M}_p(\Delta)\subseteq\mathbb{A}^{n+1}$ be the open subset consisting all non-degenerate Laurent polynomials with Newton polyhedron $\Delta$, explicitly
 \[
  \mathcal {M}_p(\Delta)= \{a\in \mathbb{A}^{n+1}\,|\, a_1\cdots a_{n+1}\prod\left(\pm 2a_1+\pm 2a_2+\cdots+\pm 2a_n+a_{n+1}\right)\ne 0\}.
 \]
 By Theorem \ref{main}(ii), we have determined $p$-adic slopes of all reciprocal zeros of $L^*(f,T)^{(-1)^n}$ for ``almost all'' $f$'s in $\mathcal {M}_p(\Delta)$ except a Zariski closed subset defined by the polynomial $h_p(\Delta)(a_1,a_2,\cdots,a_{n+1})$.
The polynomial $h_p(\Delta)$ is called a \emph{ Hasse polynomial} of the Newton polyhedron $\Delta$ . Next, we try to give an explicit Hasse polynomial. We will see that it is already quite complicated to explicitly determine the Hasse polynomial for $\Delta=\Delta(f)$ we consider, a priori for more general polyhedrons. For low dimensions, we obtain the following explicit formulae of Hasse polynomials of the Newton polyhedron $\Delta$ of the Laurent polynomials $f$ of the form~(\ref{poly}). Working on the chain level version of Dwork's trace formula, we prove
\begin{thm}
When $n=2$, a Hasse polynomial can be taken to be
 \begin{equation*}
\begin{split}
h_p(\Delta)(a_1,a_2,a_3)=\sum_{\begin{subarray}{c}0\le u+v\le {p-1\over
2}\\u,v\in \mathbb{Z}\end{subarray}}
\frac{1}{(u!v!(p-1-2u-2v)!)^2}a_1^{2v}a_2^{2u}a_3^{p-1-2u-2v}.
\end{split}
\end{equation*}
\end{thm}

For $n=3$,  it has already been a bit more complicated than the case $n=2$. Using the chain level version of Dwork's trace formula, we can easily give the condition when Newton polygon and Hodge polygon coincide at the first break point just as we treat in the case $n=2$. However, to find out when they meet at the second break point with the same method above, it needs us to compute the determinant of a matrix of size $33\times 33$ whose entries are all polynomials. In fact, it even requires a while to write down the matrix, let alone to compute the determinant.  To deal with this problem, we use the boundary decomposition theorem of Wan to divide the ``characteristic power series'' det$(I-TA_1(f))$ into ``interior pieces'' and then handle them piece by piece.

\begin{thm}\label{hass2}
When $n=3$, a Hasse polynomial can be taken to be
 \begin{equation*}
\begin{split}
h_p(\Delta)(a)=T(a) \sum_{\begin{subarray}{c}0\le u+v+w\le {p-1\over
2}\\ u,v,w\in \mathbb{Z}\end{subarray}}
\frac{1}{(u!v!w!(p-1-2u-2v-2w)!)^2}a_1^{2u}a_2^{2v}a_3^{2w}a_4^{p-1-2u-2v-2w},
\end{split}
\end{equation*}
where $a=(a_1,a_2,a_3,a_4)\in \mathcal{M}_p({\Delta})$ and $T(a)$ is explicitly presented by Formula (\ref{T}) which is essentially the determinant of a $4\times 4$ matrix whose entries are all polynomials.
\end{thm}

The rest of this paper is organized as follows. To make this paper self-contained, we first recall some baic concepts and results about $L$-functions, Newton polygon and Hodge polygon of $L$-functions. Then we review some powerful tools to study $L$-functions, such as Dwork's $p$-adic method and Wan's decomposition theory. Finally, we use these methods to investigate $L$-functions of the family of Laurent polynomials we consider in (\ref{poly}).

\section{Preliminaries}
\subsection{Exponential sums and $L$-functions}
Let $\f{q}$ be the finite field of $q$ elements with characteristic
$p$. For each positive integer $k$, let $\f{q^k}$ be the finite
extension of $\f{q}$ of degree $k$. Let $\zeta_p$ be a fixed primitive
$p$-th root of unity in the complex numbers. For any Laurent
polynomial $f(x_1,\ldots,x_n)\in
\f{q}[x_1,x_1^{-1},\ldots,x_n,x_n^{-1}]$, we form the exponential
sum
$$
S_k^*(f)=\sum_{x_1,\ldots,x_n\in
\f{q^k}^*}\zeta_p^{{\rm{Tr}}_kf(x_1,\ldots,x_n)},
$$
where $\f{q^k}^*$ denotes the multiplicative group of non-zero elements in $\f{q^k}$
and ${\rm{Tr}}_k$ denotes the trace map from $\f{q^k}$ to the prime
field $\f{p}$. To understand the sequence $S_k^*(f)\in
\mathbb{Q}(\zeta_p)$ $(1\le k<\infty)$ of algebraic integers, we
form the generating function of $S_k^*(f)$
$$
L^*(f,T)=\exp\left(\sum_{k=1}^{\infty} S_k^*(f){T^k\over k}\right),
$$
which is called the $L$-function of the exponential sum $S_k^*(f)$. The study of $L^*(f,T)$ has
fundamental importance in number theory. For example, it connects
with the zeta functions over finite fields. Consider
$$U_f(\f{q})=\{x_1,\ldots,x_n \in \f{q}^* \mid
f(x_1,\ldots,x_n)=0\}$$ the affine toric hypersurface defined by a
Laurent polynomial $$f(x_1,\ldots,x_n)\in
\f{q}[x_1,x_1^{-1},\ldots,x_n,x_n^{-1}].$$ Let $\#U_f(\f{q^k})$
denote the number of solutions of $f$ in $(\f{q^k}^*)^n$. Its zeta
function is given by
$$
Z(U_f,T)=\exp\left(\sum_{k=1}^{\infty}\#U_f(\f{q^k}) {T^k\over
k}\right).
$$
It can be easily shown that
\begin{equation}
q^k\#U_f(\f{q^k})=(q^k-1)^n+S_k^*(x_0f),
\end{equation}
and we have
$$
Z(U_f,qT)=Z(G_m^n,T)L^*(x_0f,T).
$$
Thus we see that in order to study the zeta function, it suffices to
study the $L$-function $L^*(x_0f,T)$.
Also the study of exponential sums  and the associated $L$-functions has important applications in analytic number theory, and some applied mathematics such as coding theory, cryptography, etc.

By a theorem of Dwork-Bombieri-Grothendieck, the following
generating $L$-function is a rational function:
\begin{equation}\label{rational}
L^*(f,T)=\exp\left(\sum_{k=1}^{\infty} S_k^*(f){T^k\over k}\right) =
{\prod_{i=1}^{d_1}(1-\alpha_iT)\over \prod_{j=1}^{d_2}(1-\beta_jT)},
\end{equation}
where zeros $\alpha_i\ (1\le i\le d_1)$ and poles
$\beta_j\ (1\le j\le d_2)$ are non-zero algebraic integers.
Equivalently, for each positive integer $k$, we have the formula
\begin{equation}
S_k^*(f)=\sum_{j=1}^{d_2}\beta_j^k-\sum_{i=1}^{d_1}\alpha_i^k.
\end{equation}
Thus, our fundamental question about the sums $S_k^*(f)$ is reduced
to understanding the reciprocal zeros $\alpha_i\ (1\le i\le d_1)$
and $\beta_j\ (1\le j\le d_2)$.

Without any smoothness condition of $f$, one does not even know
exactly the number $d_1$ of zeros and the number $d_2$ of poles,
although good upper bounds are available, see \cite{Bombi78}. On the
other hand, Deligne's theorem on the Riemann hypothesis \cite{Del}
gives the following general information about the nature of the
zeros and poles. For the complex absolute value $|\cdot |$, it
says
$$
|\alpha_i|=q^{u_i/2},\ |\beta_j|=q^{v_j/2},\ u_i\in
\mathbb{Z}\cap[0,2n],\ v_j\in \mathbb{Z}\cap[0,2n]
$$
where $\mathbb{Z}\cap[0,2n]$ denotes the set of integers in the
interval $[0,2n]$. Furthermore, each $\alpha_i$ (resp. each
$\beta_j$) and its Galois conjugates over $\mathbb{Q}$ have the same
complex absolute value. For each $l$-adic absolute value
$|\cdot|_l$ with prime $l\ne p$, the $\alpha_i$ and $\beta_j$ are
$l$-adic units:
$$|\alpha_i|_l=|\beta_j|_l=1.$$
For the remaining prime $p$, Deligne's integrality theorem implies that
$$
|\alpha_i|_p=q^{-r_i},\ |\beta_j|_p=q^{-s_j},\ r_i\in
\mathbb{Q}\cap[0,n],s_j\in\mathbb{Q}\cap[0,n],
$$
where the $p$-adic absolute value is normalized such that $|q|_p=1/q$.
Strictly speaking, in defining the $p$-adic absolute value, we have
tacitly chosen an embedding of the field $\bar{\mathbb{Q}}$ of
algebraic numbers into an algebraic closure of the $p$-adic number
field $\mathbb{Q}_p$. Note that each $\alpha_i$ (resp. each
$\beta_j$) and its Galois conjugates over $\mathbb{Q}$ may have
different $p$-adic absolute values. The precise version of various
types of Riemann hypothesis for the $L$-function in (\ref{rational}) is then to
determine the important arithmetic invariants $\{u_i, v_j , r_i,
s_j\}$. The integer $u_i$ (resp. $v_j$) is called the \emph{weight} of the
algebraic integer $\alpha_i$ (resp. $\beta_j$). The rational number
$r_i$ (resp. $s_j$) is called the \emph{slope} of the algebraic integer
$\alpha_i$ (resp. $\beta_j$) defined with respect to $q$. Without
any smoothness condition on $f$, not much more is known about these
weights and the slopes, since one does not even know exactly the
number $d_1$ of zeros and the number $d_2$ of poles. Under a
suitable smoothness condition, a great deal more is known about the
weights $\{u_i, v_j\}$ and the slopes $\{r_i, s_j\}$, see
Adolphson-Sperber \cite{AS89}, Denef-Loesser \cite{DL} and Wan
\cite{Wan93,Wan04}.

To investigate the slopes $\{r_i,s_j\}$, Newton polygon was introduced.

\subsection{Newton polygon and Hodge polygon}
Let $$f(x_1,\ldots,x_n)=\sum_{j=1}^{J}a_jx^{V_j},\ a_j\ne 0,$$ be a
Laurent polynomial in $\f{q}[x_1,x_1^{-1},\ldots,x_n,x_n^{-1}]$.
Each $V_j=(v_{1j},\ldots, v_{nj})$ is a lattice point in
${\mathbb{Z}}^n$ and the power $x^{V_j}$ means the product
$x_1^{v_{1j}}\cdots x_n^{v_{nj}}$. Let $\Delta(f)$ be the convex
closure in ${\mathbb{R}}^n$ generated by the origin and the lattice
points $V_j\ (1\le j\le J)$. This is called the {\it Newton
polyhedron} of $f$. If $\delta$ is a subset of $\Delta(f)$, we
define the restriction of $f$ to $\delta$ to be the Laurent
polynomial
$$
f^{\delta}=\sum_{V_j\in \delta}a_jx^{V_j}
$$
\begin{defi}
The Laurent polynomial $f$ is called {\it non-degenerate} if for
each closed face $\delta$ of $\Delta(f)$ of arbitrary dimension
which does not contain the origin, the $n$ partial derivatives
\begin{equation*}
\{ {\partial f^{\delta}\over \partial x_1},\ldots,{\partial
f^{\delta}\over \partial x_n}\}
\end{equation*}
has no common zeros with $x_1\cdots x_n\ne 0$ over the algebraic
closure of $\f{q}$.
\end{defi}

Assume now that $f$ is non-degenerate, then the $L$-function $L^*(f,
T)^{(-1)^{n-1}}$ is a polynomial of degree $n!V(f)$ by a theorem of
Adolphson-Sperber \cite{AS89} proved using $p$-adic methods, where
$V(f)$ denotes the volume of $\Delta(f)$. The complex absolute
values (or the weights) of the $n!V(f)$ zeros can be determined
explicitly by a theorem of Denef-Loeser \cite{DL} proved using
$l$-adic methods. They depend only on $\Delta$, not on the specific
$f$ and $p$ as long as $f$ is non-degenerate with
$\Delta(f)=\Delta$. Hence, the weights have no variation as $f$ and
$p$ varies. As indicated above, the $l$-adic absolute values of the
zeros are always 1 for each prime $l\ne p$. Thus, there remains the
intriguing question of determining the $p$-adic absolute values (or
the slopes) of the zeros. This is the $p$-adic Riemann hypothesis
for the $L$-function $L^*(f, T)^{(-1)^{n-1}}$. Equivalently, the
question is to determine the Newton polygon of the polynomial
\begin{equation*}
L^*(f, T)^{(-1)^{n-1}}=\sum_{i=0}^{n!V(f)}A_i(f)T^i,\ A_i(f)\in
\mathbb{Z}[\zeta_p].
\end{equation*}
The \emph{Newton polygon} of $L^*(f, T)^{(-1)^{n-1}}$, denoted by NP$(f)$, is defined to be the
lower convex closure in $\mathbb{R}^2$ of the following points
$$
(k,ord_qA_k(f)),\ k=0,1,\ldots,n!V(f).
$$
And a point in $\{(k,ord_qA_k(f))\,\mid\,\ k=1,2,\ldots,n!V(f)-1\}$ is called a \emph{break point} of the Newton polygon if the left segment and the right segment have different slopes

Let $\mathcal {N}_p(\Delta)$ be the parameter space of $f$ over
$\overline{\mathbb{F}}_p$ with fixed $\Delta(f)=\Delta$. Let
$\mathcal {M}_p(\Delta)$ be the set of non-degenerate $f$ over
$\overline{\mathbb{F}}_p$ with fixed $\Delta(f)=\Delta$. It is a
Zariski open smooth affine subset of $\mathcal {N}_p(\Delta)$. It is
non-empty if $p$ is large enough, say $p>n!V(\Delta)$. Thus $\mathcal
{M}_p(\Delta)$ is again a smooth affine variety defined over
$\f{p}$. The Grothendieck specialization theorem \cite{Wan00}
implies that as $f$ varies, the lowest Newton polygon
$$
GNP(\Delta, p)=\inf_{f\in \mathcal {M}_p(\Delta)}NP(f)
$$
exits and is attained for all $f$ in some Zariski open dense subset
of $\mathcal {M}_p(\Delta)$. The lowest polygon can then be called
the {\it generic Newton polygon}, denoted by GNP$(\Delta,p)$.

A general property is that the Newton polygon lies on or above a
certain topological or combinatorial lower bound, called the {\it
Hodge polygon} HP$(\Delta)$ which we describe bellow.

Let $\Delta$ denote the $n$-dimensional integral polyhedron
$\Delta(f)$ in $\mathbb{R}^n$ containing the origin. Let $C(\Delta)$
be the cone in $\mathbb{R}^n$ generated by $\Delta$. Then
$C(\Delta)$ is the union of all rays emanating from the origin and
passing through $\Delta$. If $c$ is a real number, we define
$c\Delta = \{cx\mid x \in \Delta\}$. For a point $u\in
\mathbb{R}^n$, the \emph{weight} $\omega(u)$ is defined to be the smallest
non-negative real number $c$ such that $u \in c\Delta$. If such $c$
does not exist, we define $\omega(u)=\infty$.

It is clear that $\omega(u)$ is finite if and only if $u \in
C(\Delta)$. If $u \in C(\Delta)$ is not the origin, the ray
emanating from the origin and passing through $u$ intersects
$\Delta$ in a face $\delta$ of codimension $1$ that does not contain
the origin. The choice of the desired $1$-codimensional face $\delta$
is in general not unique unless the intersection point is in the
interior of $\delta$. Let $\sum_{i=1}^n e_iX_i = 1$ be the equation
of the hyperplane $\delta$ in $\mathbb{R}^n$, where the coefficients
$e_i$ are uniquely determined rational numbers not all zero. Then,
by standard arguments in linear programming, one finds that the weight
function $\omega(u)$ can be computed using the formula:
\begin{equation}\label{wght}
\omega(u)=\sum_{i=1}^n e_iu_i
\end{equation}
where $(u_1,\cdots , u_n) = u$ denotes the coordinates of $u$.

Let $D(\delta)$ be the least common denominator of the rational
numbers $e_i (1 \le i \le n)$. It follows from (\ref{wght}) that for a
lattice point $u$ in $C(\delta)$, we have
\begin{equation}\label{wght1}
\omega(u) \in {1\over D(\delta)} {\mathbb{Z}}_{\ge 0},
\end{equation}where $\mathbb{Z}_{\ge 0}$ denotes the set of non-negative integers.
It is easy to show that there are lattice points $u \in C(\delta)$
such that the denominator of $\omega(u)$ is exactly $D(\delta)$.
 Let $D(\Delta)$ be the
least common multiple of all the $D(\delta)$'s:
$$
D(\Delta) = \rm{lcm}_{\delta}D(\delta),
$$
where $\delta$ runs over all the $1$-codimensional faces of $\Delta$
which do not contain the origin. Then by (\ref{wght1}), we deduce
\begin{equation}\label{wght2}
\omega({\mathbb{Z}}^n) \in {1\over D(\Delta)} {\mathbb{Z}}_{\ge
0}\cup \{\infty\}.
\end{equation}
The integer $D = D(\Delta)$ is called the \emph{denominator} of $\Delta$.
It is the smallest positive integer for which (\ref{wght2}) holds. But be
careful that there may not have a lattice point $u \in C(\Delta)$
such that the denominator of $\omega(u)$ is exactly $D(\Delta)$.

For an integer $k$, let
$$
W_{\Delta}(k) = \sharp \left\{u \in {\mathbb{Z}}^n \mid \omega(u) =
{k\over D }\right\}
$$
be the number of lattice points in ${\mathbb{Z}}^n$ with weight
$k/D$. This is a finite number for each $k$. The \emph{Hodge numbers} are defined to be
$$
H_{\Delta}(k) = \sum_{i=0}^n (-1)^i{n \choose i}W_{\Delta}(k -iD),\,k\in{\mathbb{Z}}_{\ge 0}.
$$
Hodge number $H_{\Delta}(k)$ is the number of lattice points of weight $k/D$ in a
certain fundamental domain corresponding to a basis of the $p$-adic
cohomology space used to compute the $L$-function. Thus,
$H_{\Delta}(k)$ is a non-negative integer for each $k \in
\mathbb{Z}_{\ge 0}$. Furthermore, by cohomology theory, $$H_{\Delta}(k)=0,\ for\ k > nD$$
and $$\sum_{k=0}^{nD}H_{\Delta}(k) = n!V(\Delta).$$

\begin{defi}
The {\it Hodge polygon} HP$(\Delta)$ of $\Delta$ is defined to be
the lower convex polygon in $\mathbb{R}^2$ with vertices
$$\left(\sum_{m=0}^kH_{\Delta}(k), {1\over D}\sum_{m=0}^kkH_{\Delta}(k)\right),\ k = 0, 1, 2,\ldots , nD.$$
That is, the polygon HP$(\Delta)$ is the polygon starting from the
origin and has a side of slope $k/D$ with horizontal length
$H_{\Delta}(k)$ for each integer $0 \le k \le nD$. For $k = 1, 2,\ldots , nD-1$, the point $$\left(\sum_{m=0}^kH_{\Delta}(k), {1\over D}\sum_{m=0}^kkH_{\Delta}(k)\right)$$ is called a \emph{break point} of the Hodge polygon if $H_{\Delta}(k+1)\neq 0$.
\end{defi}

The lower bound of Adolphson and Sperber \cite{AS89} says that if $f
\in \mathcal{M}_p(\Delta)$, then NP$(f) \ge $ HP$(\Delta)$ and they
have the same endpoint. The Laurent polynomial $f$ is called {\it
ordinary} if NP$(f) = $ HP$(\Delta(f))$. Combining with the definition of the generic Newton polygon, we deduce
\begin{prop}
For every prime $p$ and every $f\in \mathcal{M}_p(\Delta)$, we have
the inequalities $${\rm{NP}}(f) \ge {\rm{GNP}}(\Delta, p) \ge
{\rm{HP}}(\Delta).$$
\end{prop}

Let $\mathcal {H}_p(\Delta)$ be the moduli space of those $f\in
\mathcal{M}_p( \Delta)$ such that $\Delta(f)=\Delta$, $f$ is
nondegenerate with respect to $\Delta$ and NP$(f)=$HP$(\Delta)$. In
Dwork's terminology, $\mathcal {H}_p(\Delta)$ is called the Hasse
domain of the generic Laurent polynomial $f$ defined over
$\overline{\mathbb{F}}_p$ with $\Delta(f)$ contained in $\Delta$,
and it is a Zariski-open subset of $\mathcal{M}_p(\Delta)$ (possibly
empty). Moreover the complement of $\mathcal{H}_p(\Delta)$ in
$\mathcal {M}_p(\Delta)$ is an affine variety defined by a
polynomial in the variables $a_j$ (coefficients of $f$), called the {\it Hasse polynomial}
and denoted by $h_p(\Delta)$. Very little about Hasse polynomials is known. It is very difficult to compute Hasse polynomial in general. In
this paper, we study a family of Laurent polynomials and determine the Hasse polynomials in low dimensions.

\subsection{Diagonal local theory}
A Laurent polynomial $f$ is called {\it diagonal} if $f$ has exactly
$n$ non-constant terms and $\Delta=\Delta(f)$ is $n$-dimensional
(i.e., a simplex). In this case, the $L$-function can be computed
explicitly using Gauss sums. Let
\begin{equation}\label{diagf}
f(x)=\sum_{j=1}^na_jx^{V_j},\ a_j\in \f{q}^*,
\end{equation}
where $0,V_1,\ldots,V_n$ are the vertices of an $n$-dimensional
integral simplex $\Delta$ in $\mathbb{R}^n$. Let its vertex matrix be the non-singular $n\times n$ matrix
$$
M=(V_1,\ldots, V_n),
$$
where each $V_j$ is written as a column vector.
\begin{prop}\label{diagdeg}
For $f$ in (\ref{diagf}), $f$ is non-degenerate if and only if $p$ is
relatively prime to det$(M)$.
\end{prop}
\begin{pf}
Note that $\Delta(f)$ has only one face of dimension $n-1$ not
containing the origin. For this face, let $y_j=a_jx^{V_j}$, We have
\begin{equation}\label{diagsys}
x_i{\partial f\over \partial
x_i}=\sum_{j=1}^nV_{ji}(a_jx^{V_j})=\sum_{j=1}^nV_{ji}y_j, \ (1\le
i\le n).
\end{equation}
The $n$ partial derivatives ${\partial f\over \partial
x_1},\ldots,{\partial f\over \partial x_n}$ have no common zeros
with $x_1\cdots x_n\ne 0$ is equivalent to the $n$ linear equations
of $y_j$ ($1\le j\le n$) have no common zeros in (\ref{diagsys}), which is
equivalent to that $p$ is relatively prime to det$(M)$.

For any other face $\delta$ of dimension $m<n-1$, by a orthogonal
transformation, we can assume $\delta$ is on the hyperplane
$x_{m+1}=...=x_n=0$, which reduce to the above situation.
\end{pf}

Consider the solutions of the following linear system
\begin{equation}\label{linear}
M\left(\begin{array}{c}r_1\\\vdots \\ r_n\end{array}\right)\equiv 0\
({\rm mod }\ 1),\ r_i\ \text{rational},\ 0\le r_i<1.
\end{equation}
Let $S(\Delta)$ be the set of solutions $r$ of (\ref{linear}). It is easy to
see that  $S(\Delta)$ is a finite abelian group and its order is
precisely given by
$$|{\rm det}(M)|=n!V(\Delta).$$ Let $S_p(\Delta)$ be the prime to $p$ part of $S(\Delta)$.
It is an abelian subgroup of order equal to the prime to $p$ factor
of det$(M)$. In particular, $S_p(\Delta)=S(\Delta)$ if $p$ is
relatively prime to det$(M)$ , i.e., $f$ is non-degenerate.

Let $m$ be an integer relatively prime to the order of
$S_p(\Delta)$, then multiplication by $m$ induces an automorphism of
the finite abelian group $S_p(\Delta)$. The map is called the {\it
$m$-map} of $S_p(\Delta)$ denoted by $r\mapsto \{mr\}$, where
$$
\{mr\}=(\{mr_1\},\ldots,\{mr_n\})
$$
and $\{mr_i\}$ denotes the fractional part of the real number
$mr_i$. For each element $r\in S_p(\Delta)$, let $d(m,r)$ be the
smallest positive integer such that multiplication by $m^{d(m,r)}$
acts trivially on $r$, i.e. $$(m^{d(m,r)}-1)r\in \mathbb{Z}^n.$$ Let
$S_p(m,d)$ be the set of $r\in S(\Delta)$ such that $d(m,r)=d$, We
have the disjoint $m$-degree decomposition
$$
S_p(\Delta)=\bigcup_{d\ge 1} S_p(m,d).
$$

Let $\chi\,:\,\f{q}^*\rightarrow \mathbb{C}^*$ be a multiplicative character and let
$$
G_k(q)=-\sum_{a\in F_q^*}\chi (a)^{-k}\zeta_p^{{\rm Tr} (a)}\ (0\le
k\le q-2)
$$
be the Gauss sums. Then we have
\begin{thm}[\cite{Wan04}]\label{wan04}
\begin{equation*}
L^*(f/\f{q},T)^{(-1)^{n-1}}=\prod_{d\ge 1}\prod_{r\in
S_p(q,d)}\left(1-T^d\prod_{i=1}^n\chi
(a_i)^{r_i(q^d-1)}G_{r_i(q^d-1)}(q^d)\right)^{{1\over d}},
\end{equation*}
where $r=(r_1,\ldots, r_n)$.
\end{thm}

The Stickelberger theorem for Gauss sums is
\begin{thm}[\cite{GTM84}]\label{gtm}
Let $0 \le k \le q- 2$. Let $\sigma_p(k)$ be the sum of the
$p$-digits of $k$ in its base $p$ expansion. That is, $\sigma_p(k) =
k_0 + k_1 + k_2 + \ldots$, $k = k_0 + k_1p + k_2p^2 +\ldots $, $0
\le k_i \le p- 1$. Then, $${\rm ord}_pG_k(q) = {\sigma_p(k)\over p
-1}.$$
\end{thm}

By Theorems \ref{wan04} and \ref{gtm}, with a calculation, we have the ordinary
criterion for a diagonal Laurent polynomial $f$.
\begin{thm}[\cite{Wan04}]\label{diag}
Let $d_n(p)$ be the largest invariant factor of $S_p(\Delta)$. Let
$d_n$ be the largest invariant factor of $S(\Delta)$. If $p \equiv
1\  ({\rm mod}\ d_n(p))$, then the diagonal Laurent polynomial in
(\ref{diagf}) is ordinary at $p$. In particular, if $p \equiv 1 ({\rm mod}\
d_n)$, then the diagonal Laurent polynomial in (\ref{diagf}) is ordinary at
$p$.
\end{thm}

In order to study the (generically) ordinary property of $L$-functions and determine a Hasse polynomial $h_p(\Delta)$, we are going to briefly review Dwork's trace formula, Wan's descent theorem and Wan's decompostion theory for $L$-function.
\subsection{Dwork's trace formula}

Let $\mathbb{Q}_p$ be the field of $p$-adic numbers. Let $\Omega$
 be the completion of an algebraic closure of $\mathbb{Q}_p$. Let $q = p^a$ for some
positive integer $a$. Denote by ord the additive valuation on
$\Omega$ normalized by ord $p=1$, and denote by ord$_q$ the
additive valuation on $\Omega$ normalized by ord$_qq=1$. Let $K$
denote the unramified extension of $\mathbb{Q}_p$ in $\Omega$
 of degree $a$. Let $\Omega_1 = \mathbb{Q}_p(\zeta_p)$, where $\zeta_p$ is a primitive $p$-th root of unity. Then $\Omega_1$ is
the totally ramified extension of $\mathbb{Q}_p$ of degree $p-1$.
Let $\Omega_a$ be the compositum of $\Omega_1$ and $K$. Then
$\Omega_a$ is an unramified extension of $\Omega_1$ of degree $a$.
The residue fields of rings of algebraic integers of $\Omega_a$ and $K$ are both $\mathbb{F}_q$,
and the residue fields of rings of algebraic integers of $\Omega_1$ and $\mathbb{Q}_p$ are both
$\mathbb{F}_p$. Let $\pi$ be a fixed element in $\Omega_1$
satisfying
$$
\sum_{m=0}^{\infty}{{\pi}^{p^m}\over p^m}=0,\ {\rm ord}_p\pi
={1\over p-1}.
$$
Then, $\pi$ is a uniformizer of $\Omega_1 = \mathbb{Q}_p(\zeta_p)$
and we have
$$\Omega_1 = \mathbb{Q}_p(\pi).$$ The Frobenius
automorphism $x\mapsto x^p$ of Gal$(\f{q}/\mathbb{F}_p)$ lifts to a
generator $\tau$ of Gal$(K/\q{p})$ such that $\tau(\pi) = \pi$. If $\zeta$ is a $(q - 1)$-st root
of unity in $\Omega_a$, then $\tau (\zeta) = \zeta^p$.

Let $E(t)$ be the Artin-Hasse exponential series:
$$
E(t)=\exp(\sum_{m=0}^{\infty}{{t}^{p^m}\over p^m})=\prod_{k\ge
    1,(k,p)=1}(1-t^k)^{\mu(k)/k}
$$
where $\mu(k)$ is the M\"{o}bius function. The last product
expansion shows that the power series $E(t)$ has $p$-adic integral
coefficients. Thus, we can write
$$
E(t)=\sum_{m=0}^{\infty}\lambda_mt^m,\ \lambda_m\in \z{p}.
$$
For $0 \leq m \leq 2p-1$ (what we need below), more precise information is given by
\begin{equation}\label{lam1}
\lambda_m = {1\over m!} , \ {\rm ord}_p\lambda_m = 0, \ 0 \le m \le p
- 1.
\end{equation}
\begin{equation}\label{lam2}
\lambda_m = {1\over m!}+{1\over {p(m-p)!}}, \ {\rm ord}_p\lambda_m \geq  0,\ p \le m \le 2p
- 1.
\end{equation}
The shifted series
\begin{equation}\label{theta}
\theta(t) = E(\pi t) = \sum_{m=0}^{\infty}\lambda_m\pi^mt^m
\end{equation}
is a splitting function in Dwork's terminology. The value
$\theta(1)$ is a primitive $p$-th root of unity, which will be
identified with the $p$-th root of unit $\zeta_p$ used in our
definition of the exponential sums as given in the introduction.

For a Laurent polynomial $$f(x_1,\cdots , x_n) \in
\f{q}[x_1,x_1^{-1},\cdots , x_n, x_n^{-1} ],$$ we write $f =
\sum_{j=1}^J \bar{a}_jx^{V_j}, V_j \in \mathbb{Z}^n, \bar{a}_j \in
\f{q}$. Let $a_j$ be the Teichm\"{u}ller lifting of $\bar{a}_j$ in
$\Omega$. Thus, we have $a_j^q = a_j$. Set
\begin{equation}\label{thetaf}
F(f,x)=\prod_{j=1}^J\theta(a_jx^{V_j})
\end{equation}
\begin{equation}
F_a(f,x)=\prod_{i=0}^{a-1}F^{\tau^i}(f,x^{p^i}).
\end{equation}
Note that (\ref{theta}) implies that $F(f, x)$ and $F_a(f, x)$ are well
defined as formal Laurent series in $x_1, \cdots , x_n$ with
coefficients in $\Omega_a$.

To describe the growth conditions satisfied by $F$, write
$$
F(f,x)=\sum_{r\in \mathbb{Z}^n}F_r(f)x^r.
$$
Then from (\ref{theta}) and (\ref{thetaf}), one checks that
\begin{equation}\label{F_r}
F_r(f)=\sum_u(\prod_{j=1}^J\lambda_{u_j}a_j^{u_j})\pi^{u_1+\cdots+u_J},
\end{equation}
where the outer sum is over all solutions of the linear system
\begin{equation}\label{subs}
\sum_{j=1}^Ju_jV_j=r, \ u_j\ge 0 , \ u_j\ integral.
\end{equation}
Thus, $F_r(f)=0$ if (\ref{subs}) has no solutions. Otherwise, (\ref{F_r})
implies that
$$
{\rm ord}F_r(f)\ge {1\over p-1}\inf_{u}\{\sum_{j=1}^Ju_j\},
$$
where the inf is taken over all solutions of (\ref{F_r}).

For a given point $r \in \mathbb{R}^n$, recall that the weight
$\omega(r)$ is given by $$\omega(r) = \inf_u \left\{ \sum_{j=1}^Ju_j |
\sum_{j=1}^Ju_jV_j=r,u_j\ge 0,u_j\in \mathbb{R}\right\},$$ where the
weight $\omega(r)$ is defined to be $\infty$ if $r$ is not in the
cone generated by $\Delta$ and the origin. Thus,
\begin{equation}\label{ordF}
{\rm ord}F_r(f) \ge {\omega_(r)\over p - 1},
\end{equation}
with the obvious convention that $F_r(f) = 0 $ if $\omega(r) =
\infty$. Let $C(\Delta)$ be the closed cone generated by the origin
and $\Delta$. Let $L(\Delta)$ be the set of lattice points in
$C(\Delta)$. That is, $$L(\Delta) = \mathbb{Z}^n \cap C(\Delta).$$
For real numbers $b$ and $c$ with $0 < b \le p/(p - 1)$, define the
following two spaces of $p$-adic functions:
$$
\mathcal {L} (b,c)=\{\sum_{r\in L(\Delta)}C_rx^r\mid C_r\in
\Omega_a,\ {\rm ord}_pC_r\ge b\omega(r)+c\}
$$
$$
\mathcal{L}(b)=\bigcup_{c\in \mathbb{R}}\mathcal{L}(b,c).$$ one
checks from (\ref{ordF}) that
$$
F(f,x)\in \mathcal{L}({1\over p-1},0),\quad F_a(f,x)\in
\mathcal{L}({p\over q(p-1)},0).
$$

Define an operator $\psi$ on formal Laurent series by
$$
\psi(\sum_{r\in L(\Delta)}C_rx^r)=\sum_{r\in L(\Delta)}C_{pr}x^r.
$$
It is clear that
$$
\psi(\mathcal{L}(b,c))\subset \mathcal{L}(pb,c).
$$
It follows that the composite operator $\phi_a =\psi^a\circ F_a(f,
x)$ is an $\Omega_a$-linear endomorphism of the space
$\mathcal{L}(b)$, where $F_a(f, x)$ denotes the multiplication map
by the power series $F_a(f, x)$. Similarly, the operator $\phi_1
=\tau^{-1}\psi\circ F(f, x)$ is a semilinear ($\tau^{-1}$-linear)
endomorphsim of the space $\mathcal{L}(b)$. The operators $\psi_a^m$
and $\psi_1^m$ have well defined traces and Fredholm determinants.
The Dwork trace formula asserts that for each positive integer k,
$$
S_k^*(f)=(q^k-1)^n{\rm Tr}(\phi_a^k).
$$
In terms of $L$-function, this can be reformulated as follow.
\begin{thm}
We have
$$
L^*(f,T)^{(-1)^{n-1}}=\prod_{i=0}^n{\rm
det}(I-Tq^i\phi_a)^{(-1)^i{n\choose i}}.
$$
\end{thm}

The $L$-function is determined by the single determinant
det$(I-T\phi_a)$. For explicit calculations, we shall describe the
operator $\phi_a$ in terms of an infinite nuclear matrix. First, one
checks that $\phi_1^a=\phi_a$. We now describe the matrix form of
the operators $\phi_1$ and $\phi_a$ with respect to some orthonormal
basis. Fix a choice $\pi^{1/D}$ of $D$-th root of $\pi$ in $\Omega$.
Define a space of functions
$$
\mathcal {B}=\{\sum_{r\in L(\Delta)}C_r\pi^{\omega(r)}x^r\mid
C_r\in\Omega_a(\pi^{1/D}),\ C_r\rightarrow 0\ {\rm as}\
|r|\rightarrow \infty\}.
$$
Then, the monomials $\pi^{\omega(r)}x^r (r\in L(\Delta))$ form an
orthonormal basis of the $p$-adic Banach space $\mathcal{B}$.
Furthermore, if $b> 1/(p - 1)$, then $\mathcal{L}(b) \subseteq
\mathcal{B}$. The operator $\phi_a$ (resp. $\phi_1$) is an
$\Omega_a$-linear (resp. $\tau^{-1}$-semilinear ) nuclear
endomorphsim of the space $\mathcal{B}$. Let $\Gamma$ be the
orthonormal basis $\{\pi^{\omega(r)}x^r\}_{r\in L(\Delta)}$ of
 $\mathcal{B}$ written as a column vector. One checks that the
operator $\phi_1$ is given by
$$
\phi_1\Gamma=A_1(f)^{\tau^{-1}}\Gamma,
$$
where $A_1(f)$ is the infinite matrix whose rows are indexed by $r$
and columns are indexed by $s$. That is,
\begin{equation}\label{Af}
 A_1(f) = (a_{r,s}(f)) =
(F_{ps-r}(f)\pi^{\omega(r)-\omega(s)}).
\end{equation}

Since $\phi_a = \phi_1^a$ and $\phi_1$ is $\tau^{-1}$-linear, the
operator $\phi_a$ is given by
$$
\phi_a\Gamma=\phi_1^a\Gamma=\phi_1^{a-1}A_1^{\tau^{-1}}\Gamma=A_1^{\tau^{-a}}\cdots
A_1^{\tau^{-1}}\Gamma=A_1A_1^{\tau}\cdots A_1^{\tau^{a-1}}\Gamma.
$$
Let
$$
A_a(f)=A_1A_1^{\tau}\cdots A_1^{\tau^{a-1}}.
$$
Then, the matrix of $\phi_a$ under the basis $\Gamma$ is $A_a(f)$.
We call $A_1(f) = (a_{r,s}(f))$ the infinite semilinear Frobenius
matrix and $A_a(f)$ the infinite linear Frobenius matrix. Dwork's
trace formula can now be rewritten in terms of the matrix $A_a(f)$
as follows.
\begin{equation}
L^*(f,T)^{(-1)^{n-1}}=\prod_{i=0}^n{\rm
det}(I-Tq^iA_a(f))^{(-1)^i{n\choose i}}.
\end{equation}
We are reduced to understanding the single determinant det$(I
-TA_a(f))$.

\subsection{Nowton polygons of Fredholm determinants and a descent
theorem of Wan}

By (\ref{ordF}) and (\ref{Af}), we obtain the estimate
\begin{equation}\label{ordAf}
{\rm ord}a_{r,s}(f)\ge {\omega(ps-r)+\omega(r)-\omega(s)\over
p-1}\ge \omega(s).
\end{equation}\label{block}
Let $\xi\in \Omega$ be such that $\xi^D=\pi^{p-1}$. Then
ord$_p\xi=1/D$. By the above estimate (\ref{ordAf}), the infinite matrix $A_1(f)$ has the
block form
\begin{equation}
A_1(f)=\left( \begin{array}{ccccc}A_{00}& \xi^1A_{01}&\cdots
&\xi^iA_{0i}&\cdots \\
A_{10}& \xi^1A_{11}&\cdots
&\xi^iA_{1i}&\cdots \\
\vdots&\vdots & \ddots&\vdots & \\
A_{i0}& \xi^1A_{i1}&\cdots
&\xi^iA_{ii}&\cdots \\
\vdots&\vdots & \ddots&\vdots & \\
\end{array}\right)
\end{equation}
where the block $A_{ij}$ is a finite matrix of $W_{\Delta}(i)$ rows
and $W_{\Delta}(j)$ columns whose entries are $p$-adic integers in
$\Omega$. Now we introduce Wan's descent method to consider the chain level.

\begin{defi}
Let $P(\Delta)$ defined to be the polygon in $\mathbb{R}^2$ with
vertices (0,0) and
$$
\left( \sum_{k=0}^mW_{\Delta}(k),{1\over
D(\Delta)}\sum_{k=0}^mkW_{\Delta}(k)\right),\ m=0,1,2,\ldots.
$$
\end{defi}
The block form in (\ref{block}) and the standard determinant expansion of
the Fredholm determinant shows that we have
\begin{prop}
The Newton polygon of det$(I-TA_1(f))$ computed with respect to $p$
lies above the polygon $P(\Delta)$.
\end{prop}

Using the block form (\ref{block}) and the exterior power construction of a
semi-linear operator, one then gets the following lower bound of
Adolphson and Sperber [1] for the Newton polygon of det$(I-
TA_a(f))$.
\begin{prop}
The Newton polygon of det$(I-TA_a(f))$ computed with respect to $q$
lies above the polygon $P(\Delta)$.
\end{prop}

For the application to $L$-function, we need to use the linear
Frobenius matrix $A_a(f)$ instead of the simpler semi-linear
Frobenius matrix $A_1(f)$. In general, the Newton polygon of det$(I-
TA_a(f))$ computed with respect to $q$ is different from the Newton
polygon of det$(I - TA_1(f))$ computed with respect to $p$, even
though they have the same lower bound. Since the matrix $A_a(f)$ is
much more complicated than $A_1(f)$, especially for large $a$, we
would like to replace $A_a(f)$ by the simpler matrix $A_1(f)$. This
is not possible in general. However, if we are only interested in
the question whether the Newton polygon of det$(I - TA_a(f))$
coincides with its lower bound, the following theorem shows that we
can descend to the simpler det$(I - TA_1(f))$.

\begin{thm}[\cite{Wan04}]\label{chainlevel}
Let $\Delta(f) = \Delta$. Assume that the $L$-function $L^*(f,
T)^{(-1)^{n-1}}$ is a polynomial. Then, NP$(f)$ = HP$(\Delta)$ if
and only if the Newton polygon of det$(I - TA_1(f))$ coincides with
its lower bound $P(\Delta)$. In this case, the degree of the
polynomial $L^*(f, T)^{(-1)^{n-1}}$ is exactly $n!V(f)$.
\end{thm}

The theorem of Adolphson-Sperber shows that the polynomial condition
of Theorem \ref{chainlevel} is satisfied for every non-degenerate $f$ with $n$-dimensional $\Delta(f)$.

\subsection{Global decomposition theory}
In this subsection, we describe the basic facial decomposition theorem,
star decomposition theorem and boundary decomposition theorem from \cite{Wan93} for the Newton
polygon. This is enough to investigate the family of Laurent polynomials in
this paper, for more decomposition theorems and their applications,
we refer to \cite{Wan93,Wan04,Phong}.

\subsubsection{Facial decomposition for the Newton polygon}
Let $f(x)$ be a Laurent polynomial over $F_q$ such that $\Delta(f) =
\Delta$ is $n$-dimensional. Assume that $f$ is non-degenerate, and
$\delta_1,\ldots , \delta_h$ are all the $1$-codimensional faces of
$\Delta$ which do not contain the origin. Let $f^{\delta_i}$ be the
restriction of $f$ to the face $\delta_i$. Then,
$\Delta(f^{\delta_i} ) = \Delta_i$ is $n$-dimensional. Furthermore,
since $f$ is non-degenerate, it follows that each $f^{\delta_i}$ is
also non-degenerate by the definition of non-degenerate. Then we
have the following facial decomposition theorem.

\begin{figure}[!htbp]
\begin{center}
\includegraphics{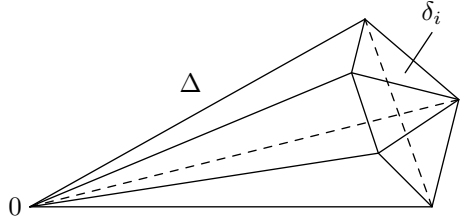}
\caption{Facial decomposition of $\Delta$}
\end{center}
\end{figure}
\begin{thm}[\cite{Wan93}] \label{facial}
Let $f$ be non-degenerate and let $\Delta(f)$ be $n$-dimensional. Then $f$ is
ordinary if and only if each $f^{\delta_i}(1\le i\le h)$ is
ordinary.
\end{thm}

\subsubsection{Star decomposition for generic Newton polygon}
By the facial decomposition, we could assume $\Delta(f)$ has only
one face $\Delta'$ of co-dimension 1 which does not contain the
origin. Let $V_1$ be an integral point on $\Delta'$. Let
$\delta_1,\ldots , \delta_h$ be the face of $\Delta'$ of
co-dimension 1 which does not contain $V_1$. For $(1\le i\le h)$, we
denote $\Delta_i$ be the convex polyhedron of $\delta_i,V_1$ and the
origin. Then we have the (closed) star decomposition:

\begin{thm}[\cite{Wan93}]\label{star} Let $f$ be
non-degenerate and let $\Delta(f)$ be $n$-dimensional with only one
face $\Delta'$ of co-dimension 1 not containing the origin. Then $f$
is generically ordinary if each $f^{\Delta_i}(1\le i\le h)$ is
generically ordinary.
\end{thm}
\begin{figure}[!htbp]
\begin{center}
\includegraphics{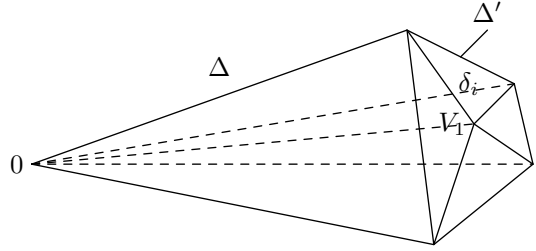}
\caption{Star decomposition of $\Delta$}
\end{center}
\end{figure}

\subsubsection{Boundary decomposition for Newton polygon}
Let $B(\Delta)$ be the unique interior decomposition of the cone $C(\Delta)$ into a union of disjoint relatively open cones. Its elements are the interiors of those closed faces in $C(\Delta)$ which contain the origin. Note that the origin itself is the unique element in $B(\Delta)$ of dimension $0$. For $\Sigma \in B(\Delta)$, denote by $A_{1}(\Sigma)$ the ``$\Sigma$'' piece $(a_{s,r}(f))$ in $A_{1}(f)$ such that $r$ and $s$ run through the cone $\Sigma$. And let $A_{1}(\Sigma,f_{\Sigma})$ be the ``interior'' piece of the Frobenius matrix $A_{1}(f_{\Sigma})$, where $f_\Sigma$ is the restriction of $f$ to the closure of $\Sigma$. The boundary decomposition theorem says
\begin{thm}[\cite{Wan93}]\label{boundary}
The Newton polygon of det$(I-tA_1(f))$ coincides with its lower bound $P(\Delta)$ if and only if the Newton polygon of det$(I-tA_{1}(\Sigma,f_{\Sigma}))$ coincides with its lower bound $P(\Delta)$ for all $\Sigma \in B(\Delta)$.
\end{thm}

\section{A family of exponential sums}

Recall the family of Laurent polynomials we consider is defined by
\begin{equation*}
f(x_1,x_2,\cdots,x_{n+1})=a_1x_{n+1}(x_1+{1\over x_1})+\cdots+a_{n}x_{n+1}(x_{n}+{1\over x_{n}})+a_{n+1}x_{n+1}+{1\over x_{n+1}}
\end{equation*}
where $a_i\in \f{q}^*,\,i=1,2,\cdots,n+1$. In this section, we study the (generically) ordinary property of $f$, give $p$-slopes of all zeros of the polynomial $L^*(f,T)^{(-1)^n}$ for generic $f$ and compute Hasse polynomials of $L^*(f,T)^{(-1)^n}$ in the low dimension cases.
\subsection{Generic Newton polygon of $f$}
For $n=1$, non-degenerate $f$ is always ordinary and the non-degenerate property is given in Proposition~\ref{ngen}. So we assume $n\ge 2$.

The Newton polyhedron of $f$, $\Delta=\Delta(f)$, has vertices $V_0=(0,\cdots,0,1)$, $V_1=(1,0,\cdots,0,1)$,
$V_2=(-1,0,\cdots,0,1)$, ..., $V_{2n-1}=(0,\cdots,0,1,1)$, $V_{2n}=(0,\cdots,0,-1,1)$, and $V_{2n+1}=-V_0$. $\Delta$ has $2^n+1$ faces of
codimension $1$. Explicitly, they are
\[
  \delta_0:\, x_{n+1}=1
\]
and
\[
 \delta_{(c_1,c_2,\cdots,c_n)}:\, 2c_1x_1+2c_2x_2+\cdots+2c_nx_n-x_{n+1}=1
\]
where $(c_1,c_2,\cdots,c_n)\in \{1,-1\}^n$. For example, vertices $V_1, V_3,\cdots, V_{2n-1}, V_{2n+1}$ determine the face $\delta_{(1,1,\cdots,1)}$,
vertices $V_2, V_3,\cdots, V_{2n-1}, V_{2n+1}$ determine the face $\delta_{(-1,1,\cdots,1)}$, and vertices $V_0, V_1, V_2,\cdots, V_{2n}$ determine the face $\delta_{0}$.

Then we have the denominator
\[
   D=D(\Delta)=1.
\]
Denote the restrictions of $f$ to these faces by
\[
 g=f^{\delta_0}=a_1x_{n+1}(x_1+{1\over x_1})+\cdots+a_{n}x_{n+1}(x_{n}+{1\over x_{n}})+a_{n+1}x_{n+1}
\]
and
\[
 f^{\delta_{(c_1,c_2,\cdots,c_n)}}=a_1x_1^{c_1}x_{n+1}+\cdots+a_nx_n^{c_n}x_{n+1}+{1\over x_{n+1}}.
\]

\begin{prop}\label{ngen}
$f$ is non-degenerate if and only if
$$
\prod_{(c_1,c_2,\cdots,c_n)\in \{\pm 1\}^n}\left(2c_1a_1+2c_2a_2+\cdots+2c_na_n+a_{n+1}\right)\ne 0.
$$
\end{prop}

\begin{pf}
We have seen that $\Delta$ has only $\delta_0, \delta_{(c_1,c_2,\cdots,c_n)} ((c_1,c_2,\cdots,c_n)\in \{\pm 1\}^n$)
faces of codimension $1$ not containing the origin. For $(c_1,c_2,\cdots,c_n)\in \{\pm 1\}^n$,
the restrictions $f^{\delta_{(c_1,c_2,\cdots,c_n)}}$ are diagonal. By Proposition~\ref{diagdeg}, a diagonal Laurent
polynomial is non-degenerate if and only if $p$ is relatively prime
to the determinent of its vertex matrix. The vertex matrix of $f^{\delta_{(c_1,c_2,\cdots,c_n)}}$ is
$$
M_{(c_1,c_2,\cdots,c_n)}=\left( \begin{array}{ccccc}
c_1 & 0 &\ldots &0 & 0\\
0 & c_2 &\ldots &0 & 0\\
\vdots & \vdots&\ddots &\vdots & \vdots\\
0 & 0 &\ldots &c_n & 0\\
1 & 1 &\ldots &1 & -1\\
\end{array} \right)
$$
which has absolute determinent $1$. Thus $f^{\delta_{(c_1,c_2,\cdots,c_n)}}$
are non-degenerate for all primes $p$.

For $g$, we have
\begin{equation}\label{equa1}
{\partial g\over \partial x_i}=a_{i}x_{n+1}(1-{1\over x^i}),\qquad i=1,2,\cdots,n
\end{equation}
\begin{equation}\label{equa2}
{\partial g\over\partial x_{n+1}}=a_1(x_1+{1\over x_1})+\cdots+a_{n}(x_{n}+{1\over x_{n}})+a_{n+1}.
\end{equation}
The system of equations (\ref{equa1}) has only non-zero solution $x_i=\pm 1, i=1,2,\cdots,n$.
Then (\ref{equa2}) gives the required condition
$$
\prod_{(c_1,c_2,\cdots,c_n)\in \{\pm 1\}^n}\left(2c_1a_1+2c_2a_2+\cdots+2c_na_n+a_{n+1}\right)\ne 0.
$$

For any other face $\delta$ of codimension larger than $1$, it must be a face of $\delta_{(c_1,c_2,\cdots,c_n)}$ for some
$(c_1,c_2,\cdots,c_n)\in \{\pm 1\}^n$. But we have seen that $f^{\delta_{(c_1,c_2,\cdots,c_n)}}$ is non-degenerate as a diagonal
Laurent polynomial. And hence all the restrictions of $f^{\delta_{(c_1,c_2,\cdots,c_n)}}$ to its faces not containing the origin are non-degenerate.
So we finish the proof of the proposition.
\end{pf}

\begin{lem}\label{eqvl}
The Laurent polynomial $f$ has the same ordinary property as its restriction $g$. That is, $f$ is ordinary if and only if $g$ is
ordinary, $f$ is generically ordinary if and only if $g$ is generically ordinary, for any prime.
\end{lem}
\begin{pf}
We have seen that $\delta_0, \delta_{(c_1,c_2,\cdots,c_n)} , (c_1,c_2,\cdots,c_n)\in \{\pm 1\}^n$, form the facial decomposition of $\Delta$. By Theorem \ref{facial}, f is (generically) ordinary if and only if $g=f^{\delta_0}$ and $f^{\delta_{(c_1,c_2,\cdots,c_n)}}$ are all (generically) ordinary. While, $f^{\delta_{(c_1,c_2,\cdots,c_n)}}$ are all diagonal whose vertex matrix have determinant $\pm 1$. So by Theorem~\ref{diag}, all $f^{\delta_{(c_1,c_2,\cdots,c_n)}}$ are ordinary for all primes $p$. And hence, the Laurent polynomial $f$ has the same ordinary property as its restriction $g$.
\end{pf}

\begin{prop}\label{gen}
For any prime, $g$ is generically ordinary. So $f$ is generically ordinary for all primes.
\end{prop}
\begin{pf}
Note that the unique $1$-codimensional face of $g$ has an interior point $V_0=(0,\cdots,0,1)$. So we use star decomposition for this face with respect to $V_0$. Then the restriction of $g$ to each part of the decomposition which is diagonal with unit vertex matrix (whose determinant is $\pm 1$) is ordinary for all primes $p$. So $g$ is generically ordinary by Theorem \ref{star}. And hence $f$ is generically ordinary by Lemma~\ref{eqvl}.
\end{pf}

In order to study the exponential sum $S^*(f)$, by Dwork's $p$-adic method, we need to investigate the slopes of reciprocal zeros and poles, i.e., Newton polygon, of its $L$-function. By Proposition \ref{gen}, the Newton polygon of generic $f$ coincides with the Hodge polygon. So our next task is to present the Hodge polygon of $\Delta=\Delta(f)$. As a part of $f$, we first compute the Hodge numbers of $g$.

\begin{prop}\label{hodgeg}
Let $\Delta'=\Delta(g)$. We have
\[
    \mathrm{Vol(\Delta')}={{2^{n}}\over {(n+1)!}}.
\]
and
\[
  W_{\Delta'}(j)=\sum_{i=0}^{n} 2^{n-i} \binom{n}{i} \binom{j}{n-i},\,j\in \mathbb{Z}_{\ge 0}.
\]
So the Hodge numbers are
\[
 H_{\Delta'}(k)=\binom{n}{k},\qquad k=0,1,\cdots,n,
\]
and
\[
H_{\Delta'}(k)=0,\qquad \mbox{for all}\,\, k\geq n+1.
\]
\end{prop}

\begin{pf}
As in the proof of generically ordinary, decomposite $\Delta'$ into $2^n$ parts (in the canonical way such that each part has vertices: the origin and $(0,\cdots,0,1)$). It is easy to see that each part has volume ${{1}\over {(n+1)!}}$. So
\[
    \mathrm{Vol(\Delta')}={{2^{n}}\over {(n+1)!}}.
\]

$W_{\Delta'}(j)$ counts the number of lattice points in the cone $C(\Delta')$ of weight $j$. From the geometry, it is easy to see that
\begin{equation*}
      \begin{array}{rl}
        W_{\Delta'}(j)= & \#\{(x_1,x_2,\cdots,x_n)\in \mathbb{Z}^n\,:\, |x_1|+|x_2|+\cdots+|x_n|\leq j\} \\
        = & \sum_{i=0}^{n} 2^{n-i} \binom{n}{i} \binom{j}{n-i} \\
      \end{array}
\end{equation*}
where the second equality follows from that we divide $x_1,x_2,\cdots,x_n$ into two parts: one part consisting of $x_i=0$, and the other part consisting
of non-zeros where each non-zero $x_i$ has $2$ choices of being negative or positive.

For $k=0,1,\cdots,n$, by the definition of Hodge numbers, we need to prove the following equality
\[
 H_{\Delta'}(k)=\sum_{i,j=0}^{n}(-1)^{i} 2^{n-j} \binom{n+1}{i} \binom{n}{j} \binom{k-i}{n-j}=\binom{n}{k}.
\]
Note that
\[
\sum_{i=0}^{n}(-1)^{i} \binom{n+1}{i}\binom{k-i}{n-j}
\]
is the coefficient of the term $x^{k-n+j}$ of the following generating function
\[
    (1-x)^{n+1}\cdot {1\over {(1-x)^{n+1-j}}}.
\]
So from
\[
    (1-x)^{n+1}\cdot {1\over {(1-x)^{n+1-j}}}=(1-x)^j,
\]
we have
\[
\sum_{i=0}^{n}(-1)^{i} \binom{n+1}{i}\binom{k-i}{n-j}=(-1)^{k-n+j}\binom{j}{n-k}.
\]
And hence,
\begin{equation*}
      \begin{array}{rl}
      H_{\Delta'}(k)= & \sum_{i,j=0}^{n}(-1)^{i} 2^{n-j} \binom{n+1}{i} \binom{n}{j} \binom{k-i}{n-j}\\
        = & \sum_{j=0}^{n}(-1)^{k+n-j} 2^{n-j} \binom{n}{j}\binom{j}{n-k} \\
       =& \binom{n}{k}\left(\sum_{j=n-k}^{n}(-1)^{k+n-j} 2^{n-j} \binom{k}{n-j}\right)\\
       =& \binom{n}{k}.
      \end{array}
\end{equation*}

From the following equality
\[
   (n+1)!\mathrm{Vol(\Delta')}=2^{n}=\sum_{k=0}^nH_{\Delta'}(k),
\]
we have determined the slopes of all the reciprocal zeros of $L^*(g,T)^{(-1)^n}$ for generic $g$.
So for any integer $k>n$, $L^*(g,T)^{(-1)^n}$ has no reciprocal zero of slope $k$, i.e., $$H_{\Delta'}(k)=0.$$
\end{pf}

\begin{rem}
  Another formula for $W_{\Delta'}(j)$ is available using the formulae given in~\cite{Zhu}
\[
  W_{\Delta'}(j)=\sum_{i=0}^{n} (-1)^{i} 2^{n-i} \binom{n}{i} \binom{n+j-i}{j},\,j=0,1,\cdots,n.
\]
\end{rem}

\begin{cor}
$W_{\Delta'}(0)=1$, $W_{\Delta'}(1)=2n+1$, $W_{\Delta'}(2)=2n^2+2n+1$, \\$W_{\Delta'}(3)={4\over 3}n^3+2n^2+{8\over 3}n+1$.
\end{cor}

\begin{cor}
\begin{itemize}
  \item[(i)] For $n=2$, we have
\begin{center}
  \begin{tabular}{|c|c|c|c|c|}
    \hline
    $j$ & 0 & 1 & 2 & 3\\
    \hline
    $W_{\Delta'}(j)$ & 1 & 5 & 13&25  \\
    \hline
    $H_{\Delta'}(j)$ & 1 & 2 &1&0 \\
    \hline
  \end{tabular}
\end{center}
  \item[(ii)] For $n=3$, we have
  \begin{center}
  \begin{tabular}{|c|c|c|c|c|c|}
    \hline
    $j$ & 0 & 1 & 2 &  3&4\\
    \hline
    $W_{\Delta'}(j)$ & 1 & 7 & 25 & 63&129 \\
    \hline
    $H_{\Delta'}(j)$ & 1 & 3 &3&1&0\\
        \hline
  \end{tabular}
\end{center}
\end{itemize}
\end{cor}

Now, we compute the Hodge numbers of $f$. So, we obtain the generic Newton polygon of $f$ and hence one main result of this paper, Theorem~\ref{main}.

\begin{thm}
Notations as above. The volume of $\Delta(f)$ is
\[
    \mathrm{Vol(\Delta)}={{2^{n+1}}\over {(n+1)!}}.
\]
And we have
\[
  W_{\Delta}(\omega)=\sum_{i=0}^n 2^{n-i}\binom{n}{i}\left(\binom{\omega}{n-i}+\binom{\omega-1}{n-i}\right),\, \omega\in \mathbb{Z}_{\ge 0}.
\]
Hence, the Hodge numbers of $f$ are
\[
 H_{\Delta}(k)=\binom{n+1}{k},\qquad k=0,1,\cdots,n,
\]
and
\[
   H_{\Delta}(k)=0,\qquad k\ge n+1.
\]
\end{thm}
\begin{pf}
$\Delta(f)$ has the same bottom as $\Delta(g)$, but twice height as $\Delta(g)$. So the volume of $\Delta(f)$ is
\[
    \mathrm{Vol(\Delta(f))}=2\mathrm{Vol(\Delta(g))}={{2^{n+1}}\over {(n+1)!}}.
\]

Separate the lattice points on the boundary of $\omega\Delta$ into two disjoint parts: one part on the hypersurface $x_{n+1}=\omega$ and the other part on the hypersurfaces $\pm 2x_{1}\pm 2x_{2}\pm\cdots\pm 2x_{n}-x_{n+1}=\omega$ with $x_{n+1}<\omega$. Translate down the second part $\omega$ units, then we get
\begin{equation*}
    \begin{array}{rl}
      W_{\Delta}(\omega)-W_{\Delta'}(\omega)= & \#\{(x_1,\cdots,x_n)\in \mathbb{Z}^n\,:\,2|x_1|+2|x_2|+\\
                                              & \qquad\cdots+2|x_n|-x_{n+1}=2\omega,\,x_{n+1}<0\} \\
      = & \#\{(x_1,\cdots,x_n)\in \mathbb{Z}^n\,:\,|x_1|+\cdots+|x_n|\le \omega-1\} \\
      = & \sum_{i=0}^n 2^{n-i}\binom{n}{i} \binom{\omega-1}{n-i}.
    \end{array}
\end{equation*}
Together with Proposition \ref{hodgeg}, we obtain
\[
  W_{\Delta}(\omega)=\sum_{i=0}^n 2^{n-i}\binom{n}{i}\left(\binom{\omega}{n-i}+\binom{\omega-1}{n-i}\right),\, \omega\in \mathbb{Z}_{\ge 0}.
\]

We notice that
\[
   W_{\Delta}(\omega)=W_{\Delta'}(\omega)+W_{\Delta'}(\omega-1)\,\, \mbox{for all }\, \, \omega\in \mathbb{Z}_{\ge 0}.
\]
So for $k=0,1,\cdots,n$,
\[
  H_{\Delta}(k)=H_{\Delta'}(k)+H_{\Delta'}(k-1)=\binom{n+1}{k}.
\]
For $k\geq n+1$, the same argument as in Proposition~\ref{hodgeg} shows
\[
 H_{\Delta}(k)=0.
\]
\end{pf}

\begin{cor}\begin{itemize}
\item[(i)] For $n=2$, we have
\begin{center}
  \begin{tabular}{|c|c|c|c|c|}
    \hline
    $j$ & 0 & 1 & 2 & 3\\
    \hline
    $W_{\Delta}(j)$ & 1 & 6 & 18&38  \\
    \hline
    $H_{\Delta}(j)$ & 1 & 3 &3&1 \\
    \hline
  \end{tabular}
\end{center}
\item[(ii)] For $n=3$, we have
\begin{center}
  \begin{tabular}{|c|c|c|c|c|c|}
    \hline
    $j$ & 0 & 1 & 2 & 3&4\\
    \hline
    $W_{\Delta}(j)$ & 1 & 8 & 32&88 &192 \\
    \hline
    $H_{\Delta}(j)$ & 1 & 4 &6&4&1 \\
    \hline
  \end{tabular}
\end{center}
\end{itemize}
\end{cor}

The generically ordinary property tells us that $f$ is ordinary on a very large subset of the space $\mathcal{M}_p(\Delta)$, a Zariski open dense subset. Next, we want to determine the explicit Zariski open subset, or equivalently to determine the explicit Hasse polynomial defining the closed complement. It is a fundamental problem in Newton-Hodge theory.

To compute the Hasse polynomial for a general Laurent polynomial, by Theorem \ref{chainlevel}, we work on
the chain level $P(\Delta)$.
We have det$(I - TA_1(f)) =
\sum_{j=0}^\infty c_jT^j$ with ${\rm{ord}}_p(c_j)\ge P(\Delta, j)$
where $P(\Delta, x)$ is the piecewise linear function describing $P(\Delta)$ by $(x,
P(\Delta, x))\in P(\Delta)$.

Recall the block form (\ref{block}) of $A_1(f)$, Consider the $k$-th vertex
of $P(\Delta)$, $( \sum_{i=0}^k W_{\Delta}(i),$\\$ 1/D \sum_{i=0}^k
iW_{\Delta}(i))$. Denote $j = \sum_{i=0}^k
W_{\Delta}(i)$, we have
$$
c_j=p^{{1\over D}\sum_{i=0}^k iW_{\Delta}(i)}{\rm det} \left(
\begin{array}{cccc}
A_{00}& A_{01}&\cdots
&A_{0k}\\
A_{10}& A_{11}&\cdots
&A_{1k} \\
\vdots&\vdots & \ddots&\vdots \\
A_{k0}& A_{k1}&\cdots&A_{kk}\\
\end{array}
\right) +\underbrace{p^{{1\over D}+{1\over D}\sum_{i=0}^k iW_{\Delta}(i)}\cdot
u_j }_{{\rm error\ term}}
$$
where $u_j$ are $p$-adic integers. Thus ord$_pc_j=P(\Delta, j)$ if
and only if
$$
h_p(\Delta,k):={\rm det} \left(
\begin{array}{cccc}
A_{00}& A_{01}&\cdots
&A_{0k}\\
A_{10}& A_{11}&\cdots
&A_{1k} \\
\vdots&\vdots & \ddots&\vdots \\
A_{k0}& A_{k1}&\cdots&A_{kk}\\
\end{array}
\right)\not \equiv 0\mod{\xi}.
$$

Hence we see that $f$ is ordinary is equivalent to $h_p(\Delta,k)\ne
0$ for each $k$. Actually we need only check $k\le n$, therefore we could take
$$h_p(\Delta)=\prod_{k=0}^nh_p(\Delta,k).$$

By Lemma \ref{eqvl}, to compute a Hasse polynomial of $\Delta(f)$, the Newton polyhedron of the Laurent polynomials we consider in this paper, is equivalent to computing a Hasse polynomial of $\Delta(g)$. So we reduce our problem to only treat
\[
  g=a_1x_{n+1}(x_1+{1\over x_1})+\cdots+a_{n}x_{n+1}(x_{n}+{1\over x_{n}})+a_{n+1}x_{n+1}.
\]

From the explicit formula of $W_{\Delta(g)}(j)$, we notice that it will be very difficult to compute the Hasse polynomials as $n$ increases. Even it is
very hard to figure out when the Newton polygon of det$(I - TA_1(g))$ with respect to $p$ coincides with $P(\Delta')$ on the $3^{rd}$ segment of slope $2$, since one has to compute the determinant of a matrix with size $(2n^2+4n+3)\times(2n^2+4n+3)$ which is already $33\times33$ in the case $n=3$.

\subsection{The case $n=2$}

Write
\begin{equation}
f(x,y,z)=az(x+{1\over x})+bz(y+{1\over y})+cz+{1\over z}
\end{equation}
where $a,b,c\in \f{q}^*$. We
have seen that $f$ is non-degenerate and generically ordinary whenever $\pm 2a\pm 2b+c\ne 0$.
So we always assume $f$ is non-degenerate, i.e. $\pm 2a\pm 2b+c\ne 0$.



By the above discuss, to compute a Hasse polynomial of $\Delta(f)$, we consider
\[
 g(x,y,z)=az(x+{1\over x})+bz(y+{1\over y})+cz
\]  Let ${\Delta'} = \Delta(g)$. The only face of
co-dimension $1$ of ${\Delta'}$ which does not contain the origin is on the plane $z=1$,
thus, we have $D'=D({\Delta'})=1$.
Also we have computed that
$W_{\Delta'}(0)=1$, $W_{\Delta'}(1)=5$, $W_{\Delta'}(2)=13$ and
$W_{\Delta'}(3)=25$, $H_{\Delta'}(0)=1$, $H_{\Delta'}(1)=2$,
$H_{\Delta'}(2)=1$ and $H_{\Delta'}(k)=0$, $\forall k\ge 3$. For
non-degenerate $g$, the end points of NP$(g)$ and
HP$(\Delta')$ coincide, i.e, $h_p(\Delta',2)=1$. Also, the
origin belongs to $\Delta(g)$, we know that NP$(g)$ passes the point
(1,0), i.e. $h_p(\Delta',0) =1$. Thus, NP$(g)=$HP$(\Delta')$ if and
only if NP$(g)$ pass the point (3, 2).
This is equivalent to
$$
h_p(\Delta',1)={\rm
det}\left(\begin{array}{cc}A_{00}&A_{01}\\A_{10}&A_{11}\end{array}\right)={\rm
det}A_{11} \not \equiv 0\mod{p}.
$$
The first equality is from the trivial calculation that $A_{00}=1$
and $A_{10}=(0,\ldots,0)^T$. $A_{11}$ is a $5\times 5
(W_{\Delta'}(1)\times W_{\Delta'}(1))$ matrix indexed by the points
of weights 1 of the form:
$$\bordermatrix{\text{}&(0,0,1)&(1,0,1)&(0,1,1) &(-1,0,1)&(0,-1,1)\cr
                (0,0,1)&\ast &  0  & 0&0&0\cr
                (1,0,1)&\ast&\ast & 0 & 0&0\cr
                (0,1,1)&\ast& 0& \ast& 0 & 0\cr
                (-1,0,1)&\ast& 0  &   0  &\ast & 0\cr
                (0,-1,1)&\ast&0&0&0&\ast}. $$
$0$'s in the matrix $A_{11}$ are because $ps-r$ of the corresponding element $a_{r,s}(g)$ is not in the cone $C(\Delta')$ (so the system (\ref{subs}) has no solution).

Recall that $A_1(g) = (a_{r,s}(g)) =
(F_{ps-r}(g)\pi^{\omega(r)-\omega(s)})$, so we have
$${\rm
det}A_{11}=p^{-5}F_{(0,0,p-1)}F_{(p-1,0,p-1)}F_{(0,p-1,p-1)}F_{(-p+1,0,p-1)}F_{(0,-p+1,p-1)}.$$
We then calculate $F_r(g)$ by (\ref{F_r}):
$$
F_{(0,0,p-1)}=p\sum_{\begin{subarray}{c}0\le u+v\le {p-1\over
2}\\u,v\in \mathbb{Z}\end{subarray}}
\lambda_v^2\lambda_{u}^2\lambda_{{p-1}-2u-2v}a^{2v}b^{2u}c^{p-1-2u-2v},
$$
\begin{equation*}
\begin{array}{ll}
F_{(p-1,0,p-1)} =p \lambda_{p-1}a^{p-1},&
F_{(0,p-1,p-1)} =p\lambda_{p-1}b^{p-1},\\
F_{(-p+1,0,p-1)}=p\lambda_{p-1}a^{p-1},&
F_{(0,-p+1,p-1)}=p\lambda_{p-1}b^{p-1},
\end{array}
\end{equation*}
where $\lambda_{m}={1\over m!}$ for any $0\leq m\leq p-1$.

Set
\begin{equation*}
\begin{split}
H(a,b,c):=\sum_{\begin{subarray}{c}0\le u+v\le {p-1\over
2}\\u,v\in \mathbb{Z}\end{subarray}}
\lambda_v^2\lambda_{u}^2\lambda_{{p-1}-2u-2v}a^{2v}b^{2u}c^{p-1-2u-2v}.
\end{split}
\end{equation*}
Here, we ignore the trivial factor $ab$ since it has already satisfied in the big space $\mathcal{M}_p({\Delta})$.
 By Lemma \ref{eqvl}, we get
\begin{thm}
A Hasse polynomial of $\Delta(f)$ (and $\Delta(g)$) can be taken as
$$
h_p(\Delta(f))(a,b,c)=h_p(\Delta(g))(a,b,c)= H(a,b,c).
$$
\end{thm}

\subsection{The case $n=3$}
In the previous subsection, we treat the simple case $n=2$. Now we deal with the case $n=3$. Consider the family of non-degenerate Laurent polynomials defined by
\begin{equation*}
f(x_1,x_2,x_3,x_{4})=a_1x_{4}(x_1+{1\over x_1})+a_{2}x_{4}(x_{2}+{1\over x_{2}})+a_{3}x_{4}(x_{3}+{1\over x_{3}})+a_{4}x_{4}+{1\over x_{4}}
\end{equation*}
where $a_i\in \f{q}^*,\,i=1,2,3,4$ and $\pm 2a_1 \pm 2a_2 \pm 2a_3 + a_4 \neq 0$.

From the discuss above, we have seen that in order to compute the Hasse polynomial of $\Delta(f)$, it is sufficient to compute the Hasse polynomial of $\Delta(g)$ where $g$ is of the following form
\[
   g(x_1,x_2,x_3,x_{4})=a_1x_{4}(x_1+{1\over x_1})+a_{2}x_{4}(x_{2}+{1\over x_{2}})+a_{3}x_{4}(x_{3}+{1\over x_{3}})+a_{4}x_{4}.
\]
And we know that NP$(g)$ and HP$(\Delta(g))$ agree at the point $(1,0)$ and the end point $(8,12)$. So we need to check when the Newton
polygon NP$(g)$ has break points $(4,3)$ and $(7,9)$. By the same reason as in the case $n=2$, it is equivalent to deciding when

$$
h_p(\Delta',1)={\rm
det}\left(\begin{array}{cc}A_{00}&A_{01}\\A_{10}&A_{11}\end{array}\right)={\rm
det}A_{11} \not \equiv 0\mod{p},
$$
and
$$
 h_p(\Delta',2)={\rm
det}\left(\begin{array}{ccc}A_{00}&A_{01}&A_{02}\\A_{10}&A_{11}&A_{12}\\A_{20}&A_{21}&A_{22}\end{array}\right)={\rm
det}\left(\begin{array}{cc}A_{11}&A_{12}\\A_{21}&A_{22}\end{array}\right) \not \equiv 0 \mod p.
$$

For $h_p(\Delta',1)$, the same argument we use in the case $n=2$ shows that
\begin{equation}\label{h1}
       h_p(\Delta',1)=\sum_{\begin{subarray}{c}0\le u+v+w\le {p-1\over
2}\\ u,v,w\in \mathbb{Z}\end{subarray}}
\lambda_u^2\lambda_{v}^2\lambda_{w}^2\lambda_{{p-1}-2u-2v-2w}a_1^{2u}a_2^{2v}a_3^{2w}a_4^{p-1-2u-2v-2w}.
\end{equation}

For $h_p(\Delta',2)$, it is already a little complicated to directly write down the matrix of size $32\times32$
\[
 \left(\begin{array}{cc}A_{11}&A_{12}\\A_{21}&A_{22}\end{array}\right).
\]
But if one does so, he will find there are a lot of 0's in the matrix. So this induces us to use the boundary
decomposition theorem \ref{boundary}.

There is only one face $\delta$ of codimension $1$ in $\Delta'$, who has vertices
$V_0=(0,0,0,1)$, $V_1=(1,0,0,1)$, $V_2=(-1,0,0,1)$, $V_3=(0,1,0,1)$, $V_4=(0,-1,0,1)$, $V_5=(0,0,1,1)$, $V_6=(0,0,-1,1)$.
Divide $\delta$ into disjoint union of relatively open faces $\delta_0$ of dimension $3$, $\delta_1,\cdots,\delta_8$ of dimension $2$,
$\delta_9,\cdots,\delta_{20}$ of dimension $1$, and vertices $\delta_{21}=V_1,\cdots,\delta_{26}=V_6$ of dimension $0$. Then the open
cones $C(\delta_i), i=1,2,\cdots,26$, and $O=(0,0,0,0)$ form an open decomposition of $C(\Delta')$.

There are only $2$ interior lattice points $V_i$ and $2V_i$ of weight $\leq 2$ in $C(\delta_i)$, for each $21\leq i\leq 26$. Taking $V_1$ as an example, we have $g=a_1x_1x_4$. One can easily figure out the matrix
$$\bordermatrix{\text{}&(1,0,0,1)&(2,0,0,2)\cr
               (1,0,0,1)&a_{11}(\delta_1,g_{\delta_1}) &  a_{12}(\delta_1,g_{\delta_1})  \cr
                (2,0,0,2)&a_{21}(\delta_1,g_{\delta_1})&a_{22}(\delta_1,g_{\delta_1})\cr}=
                \bordermatrix{\text{}&(1,0,0,1)&(2,0,0,2)\cr
               (1,0,0,1)&\lambda_{p-1}a_1^{p-1} &  \lambda_{2p-1}a_1^{2p-1}  \cr
                (2,0,0,2)&0&\lambda_{2p-2}a_1^{2p-2}\cr},$$
which has determinant
\[
      \lambda_{p-1}\lambda_{2p-2}a_1^{3p-3}.
\]
Do the same for $21\leq i\leq 26$ and by multiplying these results, we obtain the corresponding trivial factor of Hasse polynomial which is already satisfied in the space $\mathcal{M}_p({\Delta})$
\[
    a_1a_2a_3.
\]

There is only $1$ interior lattice point of weight $\leq 2$ in $C(\delta_i)$, i.e., the middle point of $\delta_i$, for each $9\leq i\leq 20$. For example, we work on $\delta_9$ the open segment with end points $(2,0,0,2)$ and $(0,2,0,2)$. So the unique interior lattice point of weight $\leq 2$ in $C(\delta_9)$ is $(1,1,0,2)$ of weight $2$ and $g_{\delta_9}=a_1x_1x_4+a_2x_2x_4$. Hence, the corresponding matrix is $1\times1$. More precisely, it is
\[
   \lambda_{p-1}^2 a_1^{p-1}a_2^{p-1}.
\]
Do the same for other $10\leq i\leq 20$ and by multiplying these results, we still obtain the corresponding trivial factor of Hasse polynomial
\[
    a_1a_2a_3.
\]

There is no interior lattice point of weight $\leq 2$ in $C(\delta_i)$, for $1\leq i\leq 8$.

There are $7$ interior lattice points $W_0=(0,0,0,2)$, $W_1=(1,0,0,2)$, $W_2=(0,1,0,2)$, $W_3=(0,0,1,2)$, $W_4=(-1,0,0,2)$, $W_5=(0,-1,0,2)$, $W_6=(0,0,-1,2)$ of weight $\leq 2$ in $C(\delta_0)$. This is the most difficult part of $h_p(\Delta',2)$.

Denote
\begin{equation*}
\begin{split}
 h_0(a_1,a_2,a_3,a_4):=\sum_{\begin{subarray}{c}0\le u+v+w\le p-1\\ u,v,w\in \mathbb{Z}\end{subarray}}&
\lambda_u^2\lambda_{v}^2\lambda_{w}^2\lambda_{2p-2-2u-2v-2w}\cdot \\& a_1^{2u}a_2^{2v}a_3^{2w}a_4^{2p-2-2u-2v-2w},\\
h_1(a_1,a_2,a_3,a_4):=\sum_{\begin{subarray}{c}0\le u+v+w\le {{p-1}\over 2}\\ u,v,w\in \mathbb{Z}\end{subarray}}&
\lambda_u\lambda_{p-1+u}\lambda_{v}^2\lambda_{w}^2\lambda_{p-1-2u-2v-2w}\cdot \\&  a_1^{p-1+2u}a_2^{2v}a_3^{2w}a_4^{p-1-2u-2v-2w},\\
h_2(a_1,a_2,a_3,a_4):=\sum_{\begin{subarray}{c}0\le u+v+w\le {{p-3}\over 2}\\ u,v,w\in \mathbb{Z}\end{subarray}}&
\lambda_u\lambda_{1+u}\lambda_{v}\lambda_{p+v}\lambda_{w}^2\lambda_{p-3-2u-2v-2w}\cdot \\&  a_1^{1+2u}a_2^{p+2v}a_3^{2w}a_4^{p-3-2u-2v-2w},\\
h_3(a_1,a_2,a_3,a_4):=\sum_{\begin{subarray}{c}0\le u+v+w\le {{p-3}\over 2}\\ u,v,w\in \mathbb{Z}\end{subarray}}&
\lambda_u\lambda_{p+1+u}\lambda_{v}^2\lambda_{w}^2\lambda_{p-3-2u-2v-2w}\cdot \\&  a_1^{p+1+2u}a_2^{2v}a_3^{2w}a_4^{p-3-2u-2v-2w},\\
h_4(a_1,a_2,a_3,a_4):=\sum_{\begin{subarray}{c}0\le u+v+w\le {{p-2}\over 2}\\ u,v,w\in \mathbb{Z}\end{subarray}}&
\lambda_u\lambda_{p+u}\lambda_{v}^2\lambda_{w}^2\lambda_{p-1-2u-2v-2w}\cdot \\&  a_1^{p+2u}a_2^{2v}a_3^{2w}a_4^{p-2-2u-2v-2w},\\
h_5(a_1,a_2,a_3,a_4):=\sum_{\begin{subarray}{c}0\le u+v+w\le {{2p-3}\over 2}\\ u,v,w\in \mathbb{Z}\end{subarray}}&
\lambda_u\lambda_{1+u}\lambda_{v}^2\lambda_{w}^2\lambda_{p-1-2u-2v-2w}\cdot \\&  a_1^{1+2u}a_2^{2v}a_3^{2w}a_4^{2p-3-2u-2v-2w},
\end{split}
\end{equation*}
and
\begin{equation*}
h_6(a_1,a_2,a_3,a_4):=h_1(a_1,a_2,a_3,a_4)-h_4(a_1,a_2,a_3,a_4),
\end{equation*}
where $\lambda_{m}={1\over m!}$ for any $0\leq m\leq p-1$ and  $\lambda_{m}={1\over m!}+{1\over {p(m-p)!}}$ for any $p\leq m\leq 2p-1$.

Computing the matrix $\left(a_{W_i,W_j}(C(\delta_0),g_{\delta_0})\right)_{0\leq i,j\leq 6}$ directly, we can notice that the matrix has a kind of period property. Then using row and column transformations, the determinant of $\left(a_{W_i,W_j}(C(\delta_0),g_{\delta_0})\right)_{0\leq i,j\leq 6}$
is reduced to the product of the determinant of the following $4\times 4$ matrix
\begin{equation*}
    \left(
      \begin{array}{cccc}
        h_0(a_1,a_2,a_3,a_4) & h_4(a_1,a_2,a_3,a_4) & h_4(a_2,a_1,a_3,a_4) & h_4(a_3,a_1,a_2,a_4) \\
        2h_5(a_1,a_2,a_3,a_4) & -h_6(a_1,a_2,a_3,a_4) & 2h_2(a_1,a_2,a_3,a_4) & 2h_2(a_2,a_3,a_1,a_4) \\
        2h_5(a_2,a_1,a_3,a_4) & 2h_2(a_2,a_1,a_3,a_4) & -h_6(a_2,a_1,a_3,a_4) & 2h_2(a_1,a_3,a_2,a_4) \\
        2h_5(a_3,a_1,a_2,a_4) & 2h_2(a_3,a_1,a_2,a_4) & 2h_2(a_3,a_2,a_1,a_4) & -h_6(a_3,a_1,a_2,a_4) \\
      \end{array}
    \right)
\end{equation*}
and
\[
 h_6(a_1,a_2,a_3,a_4)h_6(a_2,a_1,a_3,a_4)h_6(a_3,a_1,a_2,a_4).
\]
 And let
\begin{equation}\label{T}
   T(a_1,a_2,a_3,a_4)
\end{equation}
denote the result above, i.e., the determinant of $\left(a_{W_i,W_j}(C(\delta_0),g_{\delta_0})\right)_{0\leq i,j\leq 6}$.

Then, we multiple the above results and obtain a Hasse polynomial of $\Delta(g)$. And hence
\begin{thm}\label{hasse3}
A Hasse polynomial of $\Delta(f)$ can be taken as
\[
 h_{\Delta}(a_1,a_2,a_3,a_4)=h_p(\Delta',1)T(a_1,a_2,a_3,a_4),
\]
where $h_p(\Delta',1)$ is presented by Formula (\ref{h1}) and $T(a_1,a_2,a_3,a_4)$ is given by Formula (\ref{T}).
\end{thm}

Here, $T(a_1,a_2,a_3,a_4)$ in Theorem \ref{hasse3} seems still a little complicated. How to get a simple formula of $h_p(\Delta)(a_1,a_2,a_3,a_4)$, we leave it as an open problem. And for general $n\ge 4$, we ask how to give an explicit formula for the Hasse polynomial of $\Delta(f)$.

\end{document}